\documentclass[a4paper,10pt]{article}
\usepackage[english,russian]{babel}
\usepackage[cp1251]{inputenc}
\usepackage[left=44mm,right=40mm,top=60mm,bottom=52mm]{geometry}

\usepackage{amsmath}
\usepackage{amssymb}

\newtheorem{theorem}{Теорема}

\newtheorem{corollary}{Следствие}
\newtheorem{definition}{Определение}
\newtheorem{lemma}{Лемма}

\newtheorem{remark}{Замечание}

\begin{document}
\section*{Об экспоненциальной скорости сходимости распределения одной нерегенерирующей системы надёжности.}
\begin{flushright}
\textbf{Г. А. ЗВЕРКИНА}\\

{\it Российский университет транспорта (МИИТ),\\
Институт проблем управления\\
им. В.А.Трапезникова РАН\\
zverkina@gmail.com\\}
\end{flushright}
УДК 519.21\\
{\small \textbf{Ключевые слова:} восстанавливаемый элемент с тёплым резервом, скорость сходимости распределения, коэффициент надёжности, Марковские процессы, метод склеивания.
\begin{center}
\textbf{Аннотация}
\end{center}
Получены оценки сверху для скорости сходимости распределения состояния восстанавливаемого элемента с тёплым резервом в случае, когда времена работы и восстановления обоих элементов могут быть зависимы, и они ограничены сверху и снизу экспоненциальным распределением. Полученный результат позволяет оценивать скорость сходимости коэффициента готовности к стационарному значению.}

{\small \begin{center}
\textbf{Abstract}
\end{center}
{\it G. A. Zverkina. On exponential convergence rate of distribution for some non-regenerative reliability system. Fundamentalnaya i
prikladnaya matematika, vol. **, no. *, pp. ** - **}

The exponential upper bounds for the convergence rate of the distribution of restorable element with partially energized standby redundancy are founded, in the case when all working and repair times are bounded by exponential random variable (upper and lower), and working and repair times can by dependent. The bounds for the convergence rate of the availability factor are estimated.}

\section{Введение. }
\subsection{Мотивация}

Как хорошо известно, математическая теория надёжности исследует поведение систем (элементов, приборов и т.п.), которые могут отказывать в случайные моменты времени (см. например, \cite{GBS}).
Важной частью теории надёжности является исследование поведения восстанавливаемого элемента, под которым может пониматься как отдельный прибор (устройство), так и сложная техническая система.
Этот элемент безотказно работает в течение случайного времени, затем в течение случайного времени он восстанавливается и снова начинает работать.

В технических приложениях в ситуациях, когда безотказная работа исследуемого элемента очень важна, этот элемент имеет некую ``страховку'' в виде резервных элементов.

Это могут быть дополнительные элементы, ``находящиеся на складе'': они подключаются к работе {\it через некоторое время} после отказа работающего элемента.
В ряде случаев такие перерывы в работе недопустимы.
Поэтому резервные элементы могут также находиться в ``горячем'', ``тёплом'', или ``холодном'' резерве.

Если элементы находятся в ``холодном'' резерве, то они не работают, и включаются в работу только после отказа основного элемента.
Эта ситуация совпадает с поведением системы массового обслуживания с конечным числом мест для ожидания, и здесь такая задача рассматриваться не будет.

Если резервные элементы находятся в ``горячем резерве'', это означает, что они работают в полную силу (но вхолостую), и сразу после отказа основного элемента его работа переключается на один из резервных элементов -- естественно, на тот, который в момент отказа основного элемента находится в рабочем состоянии.
При этом обычно предполагается, что распределение времени безотказной работы и времени восстановления резервных элементов совпадает с соответствующими распределениями для основного элемента.
Работа или восстановление любого всех элементов системы независимы.
Понятно, такое горячее резервирование излишне затратное.
Поведение такой системы описывается поведением системы массового обслуживания с конечным количеством постоянно работающих приборов; здесь такая задача рассматриваться не будет.

Поэтому во многих технических системах резервные элементы находятся в ``тёплом резерве'', т.е. в случае исправности основного элемента они работают ``не в полную силу''.
При этом время безотказной работы резервного элемента в холостом режиме существенно увеличивается, и в случае отказа основного элемента резервный элемент очень быстро переключается в рабочий режим.
В этом случае, как правило, распределение времени безотказной работы (включая нахождение в резерве) и время восстановления резервного элемента из ``тёплого'' резерва отличается от соответствующих распределений основного элемента.

Важной характеристикой для восстанавливаемого элемента (или системы, состоящей из восстанавливаемых элементов) является {\it коэффициент надёжности} $K_\Gamma(t)$ -- вероятность того, что в момент $t$ элемент (система) исправен.

Для большинства восстанавливаемых систем существует предел $$\lim\limits _{t\to \infty}K_\Gamma(t)=K_\Gamma(\infty),$$ и важно знать, какова скорость этой сходимости.

Заметим, что $K_\Gamma(t)=\mathbb{P} \{X_t\in \mathcal{S}_0\}$, где $X_t$ -- это состояние рассматриваемой системы, а $\mathcal{S}_0$ -- подмножество из пространства состояний системы.

Поэтому для оценки скорости сходимости $K_\Gamma(t)$ к стационарному значению достаточно оценить скорость сходимости распределения сл.в. $X_t$ к стационарному распределению -- в тех случаях, когда вычисление $K_\Gamma(t)$ затруднительно или невозможно.

\section{Восстанавливаемый элемент с одним элементом тёплого резерва: случай цепи Маркова в непрерывном времени}

\subsection{Восстанавливаемый элемент без резервирования}
Для одного восстанавливаемого элемента с функцией распределения (ф.р.) времени работы $F(s)=1-e^{-\lambda s}$, и ф.р. времени восстановления $G(s)=1-e^{-\mu s}$, как известно, в случае $K_{\Gamma_1}(0)=1$ (см., например, \cite[\S 2.3]{GBS})
\begin{equation}\label{A1}
K_{\Gamma_1}(t)\stackrel{\text{\rm def}}{=\!\!\!=}P\{\mbox{в момент $t$ элемент исправен}\}=\frac{\mu + \lambda e^{-(\lambda + \mu)t}}{\mu + \lambda}.
\end{equation}

Соответственно, скорость сходимости коэффициента готовности $K_{\Gamma_1}(t)$ к его стационарному значению $K_{\Gamma_1}(\infty)\stackrel{\text{\rm def}}{=\!\!\!=}\displaystyle\frac{\mu}{\mu + \lambda}$ легко находится, она экспоненциальна:
\begin{equation}\label{sk}
K_{\Gamma_1}(t)-K_{\Gamma_1}(\infty)=\frac {\lambda e^{-(\lambda + \mu)t}}{\mu + \lambda}.
\end{equation}

\subsection{Восстанавливаемый элемент с тёплым резервом}\label{vossst}

Теперь рассмотрим восстанавливаемый элемент в ситуации, когда время работы (и пребывания в состоянии тёплого резерва) и время восстановления резервного элемента также имеют экспоненциальное распределение.

Это значит, что исправный резервный элемент имеет {\it интенсивность отказа} зависящую от состояния основного восстанавливаемого элемента.
А именно.

Первый (основной) элемент имеет ф.р. времени безотказной работы $F_1(s)=1-e^{-\lambda_1 s}$, и в случае отказа время его восстановления имеет ф.р. $G_1(s)=1-e^{-\mu_1 s}$; эти распределения не зависят от состояния резервного элемента.

То есть, если первый (основной) элемент исправен в момент времени $t$, то
\begin{equation}\label{1}
\mathbb{P}\{\mbox{основной элемент откажет за время $[t, t + \Delta]$}\}=\lambda_1\Delta + o(\Delta).\end{equation}
Аналогично, если в момент времени $t$ основной элемент восстанавливается, то
\begin{equation}\label{2}
\mathbb{P}\{\mbox{основной элемент станет исправен за время $[t, t + \Delta]$}\}=\mu_1\Delta + o(\Delta).
\end{equation}

Постоянные величины $\lambda_1$ и $\mu_1$ -- это {\it постоянные интенсивности} отказов и восстановлений соответственно.

Второй (резервный) элемент имеет интенсивность отказа, зависящую состояния основного элемента элемент:
\begin{enumerate}
\item если основной элемент исправен, то интенсивность отказа резервного элемента равна $\lambda_2$;
\item в противном случае интенсивность отказа резервного элемента равна $\Lambda$.
\end{enumerate}
Если резервный элемент неисправен, то интенсивность его восстановления равна $\mu_2$ (для простоты здесь мы полагаем, что эта интенсивность не зависит от состояния основного восстанавливаемого элемента).

Предполагается, что в начальный момент времени $t=0$ оба элемента исправны.

Состояние описываемой системы надёжности из двух элементов в момент времени $t$ обозначается парой $X_t=(i_t,j_t)$, где $i_t,j_t\in \{0;1\}$ -- первая компонента этой пары описывает состояние первого элемента, а вторая -- второго.
Исправное состояние каждого из элементов соответствует значению $0$, а неисправное состояние обозначается знаком $1$.
Процесс $X_t$ -- это цепь Маркова в непрерывном времени с пространством состояний $\mathcal{X}=\{(0,0),(1,0), (0,1), (1,1)\}$.

Коэффициент готовности этой системы -- это $K_{\Gamma_2}\stackrel{\text{\rm def}}{=\!\!\!=} \mathbb{P} \{X_t\neq (1,1)\}$.

Обозначим $p_{i,j}\stackrel{\text{\rm def}}{=\!\!\!=}\mathbb{P}\{X_t= (i,j)\}$.
Аналогично \eqref{1}--\eqref{2} и учитывая формулу полной вероятности, имеем
$$
\left\{
\begin{array}{l}
p_{0,0}(t + \Delta)=p_{0,0}(1-\lambda_1\Delta - \lambda_2 \Delta) + p_{1,0}\mu_1 \Delta + p_{0,1}\mu_2 \Delta + o(\Delta);
\\
p_{1,0}(t + \Delta)=p_{1,0}(1-\mu_1\Delta - \Lambda\Delta) + p_{0,0}\lambda_1\Delta + p_{1,1}\mu_2 \Delta + o(\Delta);
\\
p_{0,1}(t + \Delta)=p_{0,1}(1-\mu_2\Delta - \lambda_1\Delta) + p_{0,0}\lambda_2\Delta + p_{1,1}\mu_1 \Delta + o(\Delta);
\\
p_{1,1}(t + \Delta)=p_{1,1}(1-\mu_1\Delta -\mu_2\Delta) + p_{1,0}\lambda_2 \Delta + p_{0,1}\lambda_1 \Delta + o(\Delta).
\end{array}
\right.
$$
Учитывая, что $p_{0,0} + p_{1,0} + p_{0,1} + p_{1,1}=1$, получаем уравнения Колмогорова:
\begin{equation}\label{eq}
\left\{
\begin{array}{l}
\dot p_{0,0}(t)=p_{0,0}(-\lambda_1 - \lambda_2 ) + p_{1,0}\mu_1 + p_{0,1}\mu_2 ;
\\
\dot p_{1,0}(t)=p_{0,0}(\lambda_1-\mu_2) + p_{1,0}(-\mu_1 - \Lambda-\mu_2) -p_{0,1}\mu_2 + \mu_2 ;
\\
\dot p_{0,1}(t)=p_{0,0}(\lambda_2 -\mu_1) -p_{1,0} \mu_1 + p_{0,1}(-\mu_2 - \lambda_1-\mu_1) + \mu_1.
\end{array}
\right.
\end{equation}
Поскольку $p_{0,0}(t) + p_{0,1}(t)=K_{\Gamma_1}(t)$ -- коэффициент готовности основного элемента, $\alpha_1\stackrel{\text{\rm def}}{=\!\!\!=}-(\lambda_1+\mu_1)$ является собственным значением системы \eqref{eq} -- см. \eqref{A1}.

Оставшиеся два собственных значения -- это решения квадратного уравнения:
\begin{multline}\label{korni}
\alpha_{2,3}=\displaystyle\frac{-\Lambda}{2}- \displaystyle\frac{\lambda_1}{2} - \displaystyle\frac{\lambda_2 }{2}-
\displaystyle\frac{\mu_1}{2}-\mu_2\pm
\\ \\
\pm \displaystyle\frac {\sqrt{(\Lambda + \mu_1)^2 + (\lambda_1 + \lambda_2 )^2 + 2\mu_1(\lambda_1-\lambda_2 + 2\mu_2-8)- 2\Lambda(\lambda_1 + \lambda_2 )}}{2}.
\end{multline}
Заметим, что в рассмотренном случае скорость сходимости коэффициента готовности $K_{\Gamma_2}(t)=p_{0,0}(t) + p_{0,1}(t) + p_{1,0}(t)$ к стационарному значению
$$
K_{\Gamma_2}(\infty)={ \textstyle\frac{ \mu_1(\mu_1\lambda_2+\Lambda(\lambda_1+\lambda_2 )+2(\mu_2^2+\mu_1-\mu_2)) + \mu_2(\lambda_1(\lambda_1+\lambda_2 +\mu_2) +\mu_1\mu_2(2\lambda_1+\lambda_2 +\Lambda) )} {(\lambda_1 + \mu_1) (\mu_1(2 + \lambda_2 ) + \mu_2(\mu_2 + \Lambda + \lambda_1+\lambda_2 ) + \Lambda(\lambda_1 + \lambda_2))}},
$$
может быть меньше, чем в (\ref{sk}).

Действительно,
$$
K_{\Gamma_2}(t)= K_{\Gamma_2}(\infty)+ C_1 e^{\alpha_1 t}+ C_2 e^{\alpha_2 t}+ C_3 e^{\alpha_3 t},
$$
где $C_1$, $C_2$, $C_3$ зависят от начального состояния системы.
Если $\Lambda$, $\lambda_2$ и $\mu_2$ достаточно малы, то, в соответствии с \eqref{korni},
$$
\alpha_{2}\stackrel{\text{\rm def}}{=\!\!\!=}\displaystyle\frac{-\Lambda}{2}- \displaystyle\frac{\lambda_1}{2} - \displaystyle\frac{\lambda_2 }{2}-
\displaystyle\frac{\mu_1}{2}-\mu_2+ \displaystyle\frac {\sqrt{\;\;\ast\;\;\ast\;\;\ast\;\;}}{2}>\alpha_1,
$$
и скорость сходимости $ K_{\Gamma_2}(t)$ к предельному значению не менее $ C_2 e^{\alpha_2 t}$ -- естественно, в том случае, когда $C_2\neq 0$.

Напомним, что важной задачей теории надёжности является оценка скорости сходимости коэффициента готовности к стационарному состоянию не только в описанном выше простейшем (экспоненциальном) случае, но и в более сложных ситуациях.

В этой работе мы рассмотрим случай, когда распределения времени работы и восстановления обоих (основного и резервного) элементов могут зависеть друг от друга, и их распределение описывается {\it ограниченными сверху и снизу интенсивностями}.

Для оценки скорости сходимости коэффициента готовности к стационарному значению будет оценена скорость сходимости распределения {\it полного} состояния восстанавливаемой системы в метрике полной вариации.
Из оценки этой скорости сходимости автоматически вытекает оценка скорости сходимости коэффициента готовности.

Для получения оценки скорости сходимости распределения {\it полного} состояния восстанавливаемой системы в метрике полной вариации будет использован {\it метод склеивания} (coupling method) и марковизация полумарковского процесса, описывающего поведение восстанавливаемой системы с помощью пополнения пространства состояний системы.
Напомним, что вычисления будут делаться в предположении, что распределения работы и восстановления обоих элементов описываются {\it ограниченными сверху и снизу интенсивностями.}
В этом случае будет показано, что скорость сходимости параметров изучаемой восстанавливаемой системы экспоненциальна.
Такого рода результаты имеют много важных приложений -- см., например, \cite{kal1,kal2}.

\section{Восстанавливаемый элемент с тёплым резервом: полумарковский случай}
Итак, рассматривается восстанавливаемая система, состоящая из двух элементов.

Предполагается, что теперь распределения времени работы и отказа могут быть не экспоненциальными, но они являются абсолютно непрерывными.

Под {\it полным состоянием системы} мы будем понимать две пары:
\begin{equation}\label{proc}
X_t\stackrel{\text{\rm def}}{=\!\!\!=}((i_t,x_t),(j_t,y_t)),
\end{equation}
где $i_t=0$ или $i_t=1$ если первый (основной) элемент работает или нет соответственно в момент времени $t$.
А величина $x_t$ равна времени, прошедшему с последнего изменения состояния $i_t$ первого (основного) элемента до момента времени $t$.

Пара $(j_t,y_t)$ описывает состояние резервного элемента в момент времени $t$ точно таким же образом.

Заметим, что, как правило, поведение восстанавливаемой системы с резервированием описывается полумарковским (линейчатым) процессом (см., например, \cite{GK}).
Однако использование вместо функций распределения их {\it интенсивностей} позволяет рассматривать этот процесс как Марковский процесс $X_t$ на пространстве состояний $\mathcal{X}\stackrel{\text{\rm def}}{=\!\!\!=}\{0,1\} \times R_+\times \{0,1\}\times R_+$.

\subsection{Описание поведения процесса с помощью {\it интенсивностей }}
В значительном количестве публикаций по теории массового обслуживания и теории надёжности используется {\it постоянная интенсивность} поступления требований, окончания обслуживания или отказов.
Постоянная интенсивность соответствует экспоненциальному распределению (времени между поступлениями требований, времени обслуживания, времени безотказной работы или восстановления); поведение описываемых с помощью {\it постоянных интенсивностей} модели теории массового обслуживания и теории надёжности описывается цепью Маркова в непрерывном времени, как было описано в разделе \eqref{vossst}.
Тем не менее, во многих прикладных задачах возникает необходимость рассматривать произвольные распределения описанных выше величин, и в этом случае используется понятие {\it интенсивности завершения случайного периода времени} (см., например, \cite{GBS,GK}).

\begin{definition}\label{def1}
Пусть $\mathcal{T}$ -- это некоторый случайный период времени (с ф.р. $\Phi(x)$), начавшийся в момент времени $t=0$.
И пусть к моменту времени $s>0$ этот период $\mathcal{T}$ не закончился.
Тогда
$$\mathbb{P} \{\mathcal{T} \in
(s,s+\Delta]| \mathcal{T}>s\}=\displaystyle
\frac{\Phi(s+\Delta)-\Phi(s)}{1-\Phi(s)}.
$$

Если ф.р. $\Phi(x)$ абсолютно непрерывна, то
$$
\displaystyle\phi(s)\stackrel{\text{\rm def}}{=\!=} \lim\limits
_{\Delta \downarrow 0}\frac{\mathbb{P} \{\mathcal{T} \in
(s,s+\Delta]| \mathcal{T}>s\}}{\Delta} =
\frac{\Phi'(s+0)}{1-\Phi(s)}.
$$

Функция $\phi(s)$ называется { \texttt{интенсивностью}} окончания
случайного периода $\mathcal{T}$ в момент $s$ при условии, что
$\mathcal{T}\geqslant s$. \hfill \ensuremath{\triangleright}
\end{definition}

\begin{remark}\label{-intens} 
Из определения \ref{def1} следует:
$$
\mathbb{P} \{\mathcal{T}\in(s,s+\Delta]| \mathcal{T}>s\}=\phi(s)\Delta+ o(\Delta).\eqno{\triangleright}
$$
\end{remark}
\begin{remark}\label{zv9} 
Ф.р. и плотность распределения сл.в. $\mathcal{T}$ могут быть определены по интенсивности окончания периода $\mathcal{T}$:
\begin{equation}\label{zv2} 
\Phi(s)=1-\exp\left(-\int\limits_0^s \phi(u)\,\mathrm{d}\, u\right); \qquad \Phi'(s)= \phi(s)\exp\left(-\int\limits_0^s \phi(u)\,\mathrm{d}\, u\right).
\end{equation}

Легко видеть, что если ф.р. сл.в. $\mathcal{T}$ экспоненциальна, то интенсивность окончания этого периода постоянна, и наоборот.

Если ф.р. $\Phi(s)$ не является непрерывной, то также можно рассматривать интенсивность для ф.р. $\Phi(s)$, но формулы (\ref{zv2}) принимают другой вид, учитывающий скачки функции распределения. \hfill \ensuremath{\triangleright}
\end{remark}

\begin{definition}
Случайная величина $\eta$ не превосходит случайную величину $\theta$
\texttt{по распределению}, если для всех $s\in\mathbb{R}$ верно неравенство
$$
\displaystyle F_\eta(s)=\mathbb{P} \{\eta\leqslant s\}
\geqslant \mathbb{P} \{\theta\leqslant s\}=F_\theta(s) .
$$
Иначе говоря, ф.р. сл.в. $\theta$ не превосходит ф.р. сл.в. $\eta$ {\it по распределению}.
Это -- отношение порядка.
Обозначим его через $\eta\prec \theta$ -- см. \cite{sht}.\hfill \ensuremath{\triangleright}
\end{definition}
\begin{remark}\label{rem3}
Если интенсивность $\phi(s)$ окончания случайного периода $\mathcal{T}$ удовлетворяет неравенству
$$
0<c<\phi(s)<C<\infty,
$$
то 
\begin{equation}\label{sravni}
\mathcal{T}_-\prec \mathcal{T}\prec\mathcal{T}_+,
\end{equation}
где $\mathbb{P}\{\mathcal{T}_-\leqslant s\}=1-e^{-C s}$, и $\mathbb{P}\{\mathcal{T}_-\leqslant s\}=1-e^{-c s}$.
Соответственно, $\mathbb{E}\,\mathcal{T}\leqslant \displaystyle\frac{1}{c}$.

Действительно, в соответствии с \eqref{zv2}, $\mathbb{P}\{\mathcal{T}_-\leqslant s\}=1-\exp\left(-\displaystyle\int\limits_0^s C\,\mathrm{d}\,u\right)$, и $\mathbb{P}\{\mathcal{T}_-\leqslant s\}=1-\exp\left(-\displaystyle\int\limits_0^s c\,\mathrm{d}\,u\right)$, откуда и следует \eqref{sravni}.
\hfill \ensuremath{\triangleright}
\end{remark}
Наконец, после этих подготовительных сведений, опишем поведение рассматриваемой восстанавливаемой системы на языке интенсивностей.
\subsection{Предположения и обозначения}
Предполагаем, что первый (основной) элемент имеет интенсивность отказа (в работающем режиме) и интенсивность восстановления (в неработающем режиме) {\it зависящие от полного состояния системы} -- $\lambda_1(X_t)$ и $\mu_1(X_t)$ соответственно.
Точно так же, второй элемент имеет интенсивность отказа
$
\lambda_2(X_t)$, а интенсивность восстановления второго (резервного) элемента в случае отказа равна $\mu_2(X_t)$.

Предполагаем, что при восстановлении первый (основной) элемент сразу начинает работать, а резервный элемент, даже если он был исправен, переходит в состояние тёплого резерва.

Иначе говоря,
\begin{equation}\label{perex}
\left\{
\begin{array}{l}
\mathbb{P}\{i_{t+\Delta}=1, x_{t+\Delta}\in[0,\Delta), j_{t+\Delta}=j_t, y_{t+\Delta}=y_t+\Delta | i_t=0\}=
\\
\hspace{7cm}=\lambda_1(X_t)\Delta+0(\Delta);
\\
\mathbb{P}\{i_{t+\Delta}=0, x_{t+\Delta}\in[0,\Delta), j_{t+\Delta}=j_t, y_{t+\Delta}=y_t+\Delta | i_t=1\}=
\\
\hspace{7cm}=\mu_1(X_t)\Delta+0(\Delta);
\\
\mathbb{P}\{j_{t+\Delta}=0, y_{t+\Delta}\in[0,\Delta), i_{t+\Delta}=i_t, x_{t+\Delta}=x_t+\Delta | j_t=0\}=
\\
\hspace{7cm}=\lambda(X_t)\Delta+0(\Delta);
\\
\mathbb{P}\{j_{t+\Delta}=1, y_{t+\Delta}\in[0,\Delta), i_{t+\Delta}=i_t, x_{t+\Delta}=x_t+\Delta | j_t=1\}=
\\
\hspace{7cm}=\mu_2(X_t)\Delta+0(\Delta).
\end{array}
\right.
\end{equation}

Соотношения \eqref{perex} задают переходную функцию процесса $X_t$, т.е. процесс $X_t$ является Марковским процессом с пространством состояний $\mathcal{X}$ и с переходными вероятностями \eqref{perex}.

Описанная здесь модель включает в себя рассмотренный ранее экспоненциальный случай, а также ситуации, когда в случае большой продолжительности работы основного элемента резервный элемент заранее начинает работать в более интенсивном режиме, чтобы быстрее включиться в работу в случае отказа основного элемента.
Также работа основного элемента может быть менее интенсивной, если неисправен резервный элемент, чтобы исключить возможность одновременной неисправности обоих элементов, и пр..
\begin{center}
{\bf Предположения. } 
\end{center}

{\it Положим, что}
\begin{equation}\label{usl}
\begin{array}{c}
0<\lambda_m^-\leqslant \lambda_m(X_t)\leqslant {\lambda_m^+}<\infty;
\\
\\
0<{\mu_m^-}\leqslant \mu_m(X_t)\leqslant {\mu_m^+}<\infty.
\end{array}
\end{equation}
\hfill \ensuremath{\triangleright}

Введём дополнительные случайные величины $\xi_m^+$, $\xi_m^-$, $\zeta_m^+$, $\zeta_m^-$
\begin{equation}\label{oboz}
\begin{array}{ll}
\mathbb{P} \{\xi_m^+\leqslant s\}=1-e^{\int\limits _0^s -\lambda_m^+ (u) \,\mathrm{d} u}; & \mathbb{P} \{\xi_m^-\leqslant s\}=1-e^{\int\limits _0^s -\lambda_m^- (u) \,\mathrm{d} u};
\\
\mathbb{P} \{\zeta_m^+\leqslant s\}=1-e^{\int\limits _0^s -\mu_m^+ (u) \,\mathrm{d} u}; & \mathbb{P} \{\zeta_m^-\leqslant s\}=1-e^{\int\limits _0^s -\mu_m^- (u) \,\mathrm{d} u}.
\end{array}
\end{equation}
В дальнейшем мы будем рассматривать наборы таких случайных величин $\xi_m^+(k)$, $\xi_m^-(k)$, $\zeta_m^+(k)$, $\zeta_m^-(k)$, $k\in \mathbb{N}$, полагая эти сл.в. независимыми в совокупности.
\begin{remark}\label{rem4}
Обозначим $\xi_m(k)$ -- $k$-й период работы $m$-го элемента ($m=1,2$), а $\zeta_m(k)$ -- $k$-й период восстановления $m$-го элемента.

В соответствии с Замечанием \ref{rem3}, 
$$
\zeta_m(k)\prec \zeta_m^-, \qquad \xi_m(k)\prec \xi_m^-. \eqno{\triangleright}
$$
\end{remark}

\subsection{Эргодичность процесса $X_t$}
Обозначим $\mathcal{P}_t$ -- распределение процесса $X_t$ в момент времени $t$, т.е. $\mathcal{P}_t(A)=\mathbb{P}\{X_t\in A\}$, $A\in \mathcal{B}(\mathcal{X})$.
Очевидно, $\mathcal{P}_t$ зависит от $X_0$.

При этом процесс $X_t$ не является {\it регенерирующим}, и стандартные методы доказательства его эргодичности как регенерирующего процесса неприменимы (см. \cite{zv1,zv2}).

\begin{lemma}\label{lemma}
В условиях \eqref{usl} процесс $X_t$ является эргодическим, т.е. существует инвариантная вероятностная мера $\mathcal{P}$ на пространстве $\mathcal{X}$ такая, что $\mathcal{P}_t\Longrightarrow\mathcal{P}$ для всех начальных состояний $X_0$ процесса $X_t$.
\end{lemma}
\textbf{Доказательство. }
Обозначим (для $m=1,2,\ldots$)
$$
\begin{array}{cc}
\theta_1 \stackrel{\text{\rm def}}{=\!=} \inf\{t>0:\; i_t=0\}, & \theta_1' \stackrel{\text{\rm def}}{=\!=} \inf\{t>\theta_1:\; i_t=1\}, \\
\theta_{m+1} \stackrel{\text{\rm def}}{=\!=} \inf\{t>\theta_m':\; i_t=0\}, & \theta_{m+1}' \stackrel{\text{\rm def}}{=\!=} \inf\{t>\theta_{m+1}':\; i_t=0\}.
\end{array}
$$

$\theta_k$ -- это момент $k$-го восстановления основного элемента.

Из замечаний \ref{rem3} и \ref{rem4} получаем, что
\begin{equation}\label{sr}
\mathbb{E}\,\theta_1\leqslant \displaystyle\frac{1}{\lambda^-_1}+\displaystyle\frac{1}{\mu^-_1};\qquad \mathbb{E}\,(\theta_{k+1}-\theta_k)\leqslant \displaystyle\frac{1}{\lambda^-_1}+ \displaystyle\frac{1}{\mu^-_1}.
\end{equation}

Возьмём достаточно малое $\varepsilon>0$ и обозначим $$\mathcal{S}^\varepsilon\stackrel{\text{\rm def}}{=\!\!\!=} \{(1,x,1,y):\, \max\{x,y\}<\varepsilon\}; \qquad\mathcal{S}^\varepsilon\subset \mathcal{X}.
$$

В момент времени $\theta_k$ ($k>1$) первый (основной) восстанавливаемый элемент переходит в рабочее состояние; $i_{\theta_k}=0$, $x_{\theta_k}=0$.

Относительно второго (резервного) элемента возможны два случая.

1. Если $j_{\theta_k}=1$, т.е. второй элемент неисправен, то остаточное время его пребывания в неисправном состоянии {\it по распределению} меньше сл.в. $\zeta_2^-$, и с вероятностью большей, чем $\varpi_1 \stackrel{\text{\rm def}}{=\!=} 1-e^{-\varepsilon\mu^-_2}$, второй элемент перейдёт в исправное состояние за время, меньшее $\varepsilon$, т.е. до момента времени $\theta_k+\varepsilon$ процесс $X_t$ окажется в множестве $\mathcal{S}^\varepsilon$.

2. Если $j_{\theta_k}=0$, т.е. второй элемент исправен, то его остаточное время пребывания в исправном состоянии {\it по распределению} меньше сл.в. $\xi^-$, а последующее время восстановления {\it по распределению} меньше сл.в. $\zeta_2^-$.
Поэтому с вероятностью большей, чем $\varpi_2 \stackrel{\text{\rm def}}{=\!=} \left(1-e^{-\frac{\varepsilon}{2}\mu^-_2}\right)\left(1-e^{-\frac{\varepsilon}{2}\lambda^-_2 }\right)$, второй элемент откажет быстрее, чем через время $\frac\varepsilon2$, и затем восстановится быстрее, чем за время $\frac\varepsilon2$.
То есть с вероятностью большей, чем $\varpi_2$, к моменту времени $\theta_k+\varepsilon$ процесс $X_t$ окажется в множестве $\mathcal{S}^\varepsilon$.

Т.е. с вероятностью большей, чем $\varpi \stackrel{\text{\rm def}}{=\!=} \min\{\varpi_1,\varpi_2\}>0$, при некотором $\vartheta_k\in [\theta_k,\theta_k+\varepsilon]$ процесс $X_{\vartheta_k}$ окажется в маленьком подмножестве $\mathcal{S}^\varepsilon$ пространства состояний $\mathcal{X}$.

Значит, время между последовательными попаданиями процесса $X_t$ в множество $\mathcal{S}^\varepsilon$ является геометрической суммой сл.в. с конечными средними (см. \eqref{sr}), и имеет конечное математическое ожидание.
Поэтому эргодичность процесса $X_t$ следует из принципа Харриса-Хасьминского.
\hfill \ensuremath{\blacksquare}
\section{Основной результат}
\begin{theorem}\label{th}
Если выполнены условия \eqref{usl}, то для любого начального состояния $X_0$ процесса $X_t$, {\it можно вычислить} числа $\alpha>0$ и $\mathfrak{K}=\mathfrak{K}(\alpha, \lambda^-_m, \lambda^+_m,\linebreak \mu^-_m,\mu^+_m)$ такие, что
$
\|\mathcal{P}_t-\mathcal{P}\|_{TV}\leqslant \mathfrak{K}e^{-\alpha t}
$
для всех $t\geqslant 0$.

Здесь $\|\cdot\|_{TV}$ обозначает метрику полной вариации.
\end{theorem}
\begin{corollary}
В условиях Теоремы \ref{th} верно неравенство
$$|K_{\Gamma_2}(t)-K_{\Gamma_2}(\infty)|\leqslant \mathfrak{K}e^{-\alpha t},$$
поскольку коэффициент готовности рассматриваемой восстанавливаемой системы равен
$$
K_{\Gamma_2}(t)=\mathbb{P}\{X_t\notin\{((1,\cdot),(1,\cdot))\}\}= 1-\mathcal{P}_t((1,\cdot),(1,\cdot)).\eqno\triangleright
$$
\end{corollary}

\noindent{\bf Доказательство Теоремы \ref{th}} -- это построение алгоритма вычисления чисел $\alpha>0$ и $\mathfrak{K}$, и этот алгоритм основан на {\it методе склеивания} Марковских процессов -- see, e.g., \cite{zv6,zv14}.

В этой работе не будет дано полное описание этого алгоритма, поскольку оно занимает много места и требует большого количества вычислений. Мы изложим только основные идеи доказательства.

Напомним, что ``классический'' метод склеивания неприменим для исследуемой модели, поскольку её поведение описывается Марковским процессом \emph{в непрерывном времени}.
Поэтому будет использовано понятие успешной склейки, предложенное в \cite{zv5}.

Также будет использовано понятие ``{\it общей части распределений}'', т.е. $\kappa \stackrel{\text{\rm def}}{=\!=} \displaystyle\int\limits_\mathbb{R}\min\{f_1(u), f_2(u)\}\,\mathrm{d}u$, где $f_m(\cdot)$ являются плотностями распределений.
В доказательстве применяется \emph{основная лемма склеивания} в её простейшем виде (см., например, \cite{zv13,verbut}).
О применении метода склеивания в теории массового обслуживания см., например, \cite{zv7,GZ}).

\subparagraph{1. Успешная склейка.}
Пусть $X_t$ и $\widehat{X}_t$ -- два независимых Марковских процесса с одной и той же переходной функцией \eqref{perex}, но с разными начальными состояниями в момент времени $t=0$.

Пусть на некотором вероятностном пространстве построены (зависимые) процессы $Y_t=((i_t,x_t),(j_t,y_t))$ и $\widehat{Y}_t=((\widehat{\imath}_t,\widehat{x}_t),(\widehat{\jmath}_t,\widehat{y}_t))$ таким образом, что:
\begin{enumerate}
\item $Y_t\stackrel{\mathcal{D}}{=}X_t$ и $\widehat{Y}_t\stackrel{\mathcal{D}}{=}\widehat{X}_t$ для всех \emph{фиксированных} $t$;
\item $\mathbb{P} \{\tau(X_0,\widehat{X}_0)<\infty\}=1$, где $\tau(X_0,\widehat{X}_0)=\tau(Y_0,\widehat{Y}_0)=\inf\{t>0:\; Y_t=\widehat{Y}_t\}$.
\end{enumerate}

Такая пара процессов $Y_t=((i_t,x_t),(j_t,y_t))$ и $\widehat{Y}_t=((\widehat{\imath}_t,\widehat{x}_t),(\widehat{\jmath}_t,\widehat{y}_t))$ называется \emph{успешной склейкой} процессов $X_t$ и $\widehat{X}_t$, а $\tau(X_0,\widehat{X}_0)$ называется\emph{ моментом склеивания}.

Нашей задачей будет конструирование успешной склейки и оценка экспоненциальных моментов случайной величины $\tau(X_0,\widehat{X}_0)$.

\subparagraph{2. Конструирование процессов $Y_t$ и $\widehat{Y}_t$.}
Для конструирования процессов $Y_t$ и $\widehat{Y}_t$ возьмём пространство $\prod\limits_{k=0}^\infty (\Omega_k,\mathbb{P} _k,\mathcal{B}_k)$, где $(\Omega_k,\mathbb{P} _k,\mathcal{B}_k)$ -- некоторое вероятностное пространство, на котором в дальнейшем будут строиться нужные для конструкции успешной склейки случайные величины.

Конструирование будет проводиться пошагово, в те моменты времени $t_1$, $t_2$, $\ldots$, когда меняется одна из компонент $i,j,\widehat{\imath},\widehat{\jmath}$ пары $(Y_t,\widehat{Y}_t)$.

Пусть в момент $t_k$ известно состояние обоих процессов и, следовательно, известны интенсивности окончания пребывания обоих приборов в том состоянии, в котором они находятся.
То есть известны совместные распределения этих остаточных времён пребывания приборов в их состояниях. 
То есть, в соответствии с определением состояния процесса \eqref{proc}, в момент $t_k$ известны величины $i_{t_k},j_{t_k},\widehat{\imath}_{t_k},\widehat{\jmath}_{t_k}$ и совместное распределение остаточных времён сл.в. $x_{t_k},y_{t_k},\widehat{\imath}_{t_k},\widehat{\jmath}_{t_k}$.
Это -- четыре случайные величины, и известно их совместное распределение.

На этом шаге используется одно из пространств $(\Omega_k,\mathbb{P} _k,\mathcal{B}_k)$, и на нём строятся эти четыре случайные величины -- таким образом, что их совместное распределение совпадает с заданным.
После этого выбирается минимальная из этих четырёх случайных величин, пусть её значение есть $\varsigma_k$.

В момент времени $t_{k+1}\stackrel{\text{\rm def}}{=\!\!\!=} t_k+\varsigma_k$ происходит следующее изменение одной из компонент $i,j,\widehat{\imath},\widehat{\jmath}$.

В момент $t_{k+1}$ процедура повторяется.

Понятно, что предложенная здесь процедура конструирования случайных процессов $Y_t$ и $\widehat{Y}_t$ не гарантирует совпадения в какой-либо момент времени этих процессов.
\subparagraph{3. Основная лемма склеивания (см., например, \cite{zv13,verbut}).}
Здесь мы даём Основную лемму склеивания в наипростейшей формулировке (см., например, \cite{GZ1}).

\begin{lemma}\label{osn}
Если сл.в. $\vartheta_1$ и $\vartheta_2$ имеют ф.р. $\Phi_1(s)$ и $\Phi_2(s)$ соответственно, и их общая часть $\kappa\stackrel{\text{\rm def}}{=\!\!\!=} \int\limits_\mathbb{R} \min\{\Phi_1'(s), \Phi_2'(s)\}\,\mathrm{d}\, s>0$, то на некотором вероятностном пространстве можно построить сл.в. $\widehat{\vartheta}_1$ и $\widehat{\vartheta}_2$ таким образом, что
\begin{enumerate}
\item $\widehat{\vartheta}_1\stackrel{\mathcal{D}}{=} {\vartheta}_1$, $\widehat{\vartheta}_2\stackrel{\mathcal{D}}{=} {\vartheta}_2$;
\item $\mathbb{P} \{\widehat{\vartheta}_1= \widehat{\vartheta}_2\}=\kappa$. \hfill \ensuremath{\triangleright}
\end{enumerate}
\end{lemma}
\begin{remark}\label{lem+}
Утверждение Леммы \ref{osn} естественным образом переносится на любое конечное количество случайных величин, а именно.

Пусть $\vartheta_1$, $\vartheta_2$, $\ldots$, $\vartheta_n$ -- сл.в. с плотностями распределения $\varphi_1(s)$, $\varphi_2(s)$, $\ldots$, $\varphi_n(s)$ соответственно, и $\kappa \stackrel{\text{\rm def}}{=\!\!\!=}\int\limits_\mathbb{R}\min\{\varphi_1(s), \varphi_2(s), \ldots, \varphi_n(s)\}>0$.
Тогда на некотором вероятностном пространстве можно построить сл.в. $\widehat{\vartheta}_1(s)$, $\widehat{\vartheta}_2(s)$, $\ldots$, $\widehat{\vartheta}_n(s)$ такие, что
\begin{enumerate}
\item $\widehat{\vartheta}_i\stackrel{\mathcal{D}}{=} {\vartheta}_i$, $i=1,2,\ldots n$;
\item $\mathbb{P} \{\widehat{\vartheta}_1= \widehat{\vartheta}_2=\ldots=\widehat{\vartheta}_n\}=\kappa$.
{ Доказательство} этого факта повторяет доказательство Леммы \ref{osn} из \cite{GZ1}.\hfill \ensuremath{\triangleright}
\end{enumerate}
\end{remark}

\subparagraph{4. Склеивание процессов $Y_t$ и $\widehat{Y}_t$.}

Сначала рассмотрим моменты $\theta_k$ -- моменты восстановления основного элемента процесса $Y_t$; $\theta_k\prec \sum\limits_{i=1}^k(\xi_1^-(i)+\zeta_1^-(i))$ (см. Лемму \ref{lemma}).

Аналогично тому, как это было показано в Лемме \ref{lemma}, для некоторого фиксированного $\varepsilon>0$ можно утверждать следующее:

\begin{enumerate}
\item С вероятностью большей, чем $$\pi_1\stackrel{\text{\rm def}}{=\!\!\!=} (1-e^{-\frac\varepsilon2\lambda_2^-})(1-e^{-\frac\varepsilon2\mu_2^-}) $$ на интервале $[\theta_k,\theta_k+\varepsilon]$ процесс $Y_t$ окажется в множестве $\mathcal{S}^\varepsilon$.
\item С вероятностью большей, чем $$\pi_2\stackrel{\text{\rm def}}{=\!\!\!=} ( 1-e^{-\frac\varepsilon2 \lambda_1^-})( 1-e^{-\frac\varepsilon2 \mu_1^-})( 1-e^{-\frac\varepsilon2 \Lambda^-})( 1-e^{-\frac\varepsilon2 \lambda_2^-})( 1-e^{-\frac\varepsilon2 \mu_2^-})$$
на интервале $[\theta_k,\theta_k+\varepsilon]$ процесс $\widehat{Y}_t$ окажется в множестве $\mathcal{S}^\varepsilon$.
\end{enumerate}
То есть с вероятностью большей, чем $\pi_1\pi_2$, на интервале $[\theta_k,\theta_k+\varepsilon]$ оба процесса ${Y}_t$ и $\widehat{Y}_t$ окажутся в множестве $\mathcal{S}^\varepsilon$.

В этом множестве $\mathcal{S}^\varepsilon$ все приборы исправны, и остаточное время пребывания всех приборов в исправном состоянии имеет общую часть 
\begin{multline*}
\kappa_1\stackrel{\text{\rm def}}{=\!\!\!=}\int \limits_0^\infty \min \{\lambda_1(X_{\theta_k+\varepsilon+s}) e^{-\lambda_1(X_{\theta_k+\varepsilon+s})}, \lambda_2(X_{\theta_k+\varepsilon+s})e^{-\lambda_2(X_{\theta_k+\varepsilon+s})} \}\,\mathrm{d}\, s \geqslant
\\
\geqslant\frac{\min\{\lambda_1^-,\lambda_2^-\}} {\max\{\lambda_1^+,\lambda_2^+\}}>0.
\end{multline*}

Поэтому на одном вероятностном пространстве можно построить четыре сл.в. -- остаточные времена пребывания всех элементов в исправном состоянии -- таким образом, что они совпадут с вероятностью не меньшей, чем $\kappa_1$.
То есть оба процесса окажутся одновременно в состоянии $(0,0,0,0)$ через время $\eta$, которое по распределению меньше, чем сл.в. с ф.р. $\Phi(s)=1-e^{-\min\{\lambda_1^-,\lambda_2^-,\mu_1^-,\mu_2^-\}s}$.

На самом деле можно сначала ``склеить'' состояния первых (основных) элементов обоих процессов, а затем -- состояния резервных элементов, но описание этой процедуры существенно более сложное.

Также можно рассмотреть моменты $\theta_k'$ -- моменты отказов первого (основного) элемента.

Опять повторяя рассуждения Леммы \ref{lemma}, можно оценить вероятность того, что на интервале $[\theta_k',\theta_k'+\varepsilon]$ все элементы обоих систем окажутся неисправны, и снова можно применить модифицированную основную Лемму склеивания \ref{osn} и показать, что с вероятностью большей, чем некоторое $\kappa_2>0$, процессы совпадут.

Предложенные здесь оценки вероятностей попадания пары процессов в некоторое множество не оптимальны, они могут быть улучшены с помощью рассмотрения различных комбинаций состояний элементов обоих процессов и различных способов выбора величины $\varepsilon$ и долей этой величины вместо $\frac{\varepsilon}{2}$.

\subparagraph{Выводы.} Итак, с некоторой вероятностью $\tilde{\kappa}$, оценка которой была дана здесь достаточно грубо, можно утверждать, что на любом интервале $[\theta_k,\theta_{k+1}]$ процессы $Y_t$ и $\widehat{Y}_t$ совпадут.
Момент совпадения $\tau=\tau(X_0,\widehat{X}_0)$.

Поэтому можно применить основное неравенство склеивания для построенной успешной склейки:
\begin{multline}\label{zv12} 
|\mathbb{P} \{X_t\in S\} - \mathbb{P} \{\widehat{X}_t\in S\}|= |\mathbb{P} \{Y_t\in S\} - \mathbb{P} \{\widehat{Y}_t\in S\}|=
\\
=|\mathbb{P} \{Y_t\in S\;\&\; \tau>t\} - \mathbb{P} \{\widehat{Y}_t\in S\;\&\; \tau>t\} +
\\
+ |\mathbb{P} \{Y_t\in S\;\&\; \tau<t\} - \mathbb{P} \{\widehat{Y}_t\in S\;\&\; \tau<t\}|=
\\
=|\mathbb{P} \{Y_t\in S\;\&\; \tau>t\} - \mathbb{P} \{\widehat{Y}_t\in S\;\&\; \tau>t\} |\leqslant \mathbb{P} \{\tau>t\}\leqslant
\\
\leqslant \mathbb{P} \{e^{\alpha\tau}>e^{\alpha t}\}\leqslant \frac{\mathbb{E} \,e^{\alpha\tau}}{e^{\alpha t}},
\end{multline}
поскольку при $\tau<t $ ``склеивание'' процессов $Y_t$ и $\widehat{Y}_t$ уже произошло, их распределения в момент $t$ одинаковы, и $\mathbb{P} \{Y_t\in S \,\&\, \tau<t\}=\mathbb{P} \{\widehat{Y}_t\in S\,\&\, \tau<t\}$; здесь $S\in\mathcal{B} (\mathcal{X} )$.

Таким образом, осталось оценить $\mathbb{E} \,e^{\beta\tau}$ и выяснить, при каких $\beta>0$ это математическое ожидание конечно.

Учитывая, что имеется равномерная (не зависящая от начальных условий) оценка
$$
\tau(X_0,\widehat{X}_0)\prec\theta_1+\sum\limits _{k=1}^\nu (\theta_{k+1}-\theta_k)+\eta\prec\eta+\sum\limits_{k=1}^\nu(\xi_k^-+\zeta_k^-),
$$
где $\theta_{k+1}-\theta_k\prec \xi_1^-+\zeta_1^- \stackrel{\mathcal{D}}{=} 
\xi_k^-+\zeta_k^-$, а также $\mathbb{P} \{\nu>n\}\leqslant (1-\tilde{\kappa}) ^n$, можно оценить величину $\alpha$ и $\mathfrak{K}=\mathbb{E} \,e^{\alpha\tau}$.

Поскольку оценка \eqref{zv12} равномерна по начальным условиям процессов $X_t$ и $\widehat{X}_t$, она верна и в том случае, когда процесс $\widehat{X}_t$ имеет в качестве начального распределения стационарное распределение $\mathcal{P}$, откуда следует
$$
\sup\limits_{S\in \mathcal{B}(\mathcal{X})}|\mathbb{P} \{X_t\in S\}-\mathcal{P}(S)|=\|\mathcal{P}_t(S)-\mathcal{P}(S)\|_{TV}\leqslant \frac{\mathbb{E} \,e^{\alpha\tau}}{e^{\alpha t}}=\mathfrak{K}e^{-\alpha t}.
$$

\section{Заключение}
Представленная здесь конструкция успешной склейки существенно упрощена и оценки, которые можно получить по описанной схеме, достаточно грубые.
Однако использование предложенного подхода с анализом всех возможных ситуаций, подходящих для склеивания, а также возможность учесть некоторые специфические характеристики изучаемой восстанавливаемой системы и начальное состояние процесса $X_t$ позволяют существенно улучшить предложенную здесь оценку.

\paragraph{Благодарности.} Автор бллагодарит Э.Ю.~Калимулину, В.В.~Козлова и А.Ю.Веретенникова за ценные рекомендации при подготовке статьи. 
Работа выполнена при финансовой поддержке РФФИ (проект № 17-01-00633 А).

\end{document}


\section*{On exponential convergence rate of distribution for some non-regenerative reliability system.}
\begin{flushright}
\textbf{G. A. ZVERKINA}\\

{\it Moscow State University of \\Railway Engineering (MIIT),\\
V.A. Trapeznikov Institute of Control\\
Sciences of the Russian\\
Academy of Sciences\\
zverkina@gmail.com\\}
\end{flushright}
UDC 519.21\\
{\small \textbf{Key words :} restorable element with a warm reserve, convergence rate of distribution, availability factor, Markov processes, copling method.
\begin{center}
\textbf{Abstract}
\end{center}
The exponential upper bounds for the convergence rate of the distribution of restorable element with partially energized standby redundancy are founded, in the case when all working and repair times are bounded by exponential random variable (upper and lower), and working and repair times can by dependent. The bounds for the convergence rate of the availability factor are estimated.}

\section{Introduction. }
\subsection{Motivation}

It is well known, that the mathematical theory of reliability studies the behavior of systems (elements, devices, etc.) which can fail at random times (see, e.g., \cite{GBS}).
The study of the behavior of the recoverable element is an important part of the reliability theory.
This element can be either a separate device or a complex technical system.
This item works flawlessly for a random time, then it recovers and starts to work again for a random time.

In technical applications, in situations where the failure of the studied element is very important, this element has a certain ``insurance'' in the form of reserve elements.

These may be additional items that are ``in stock'': they starts to work {\it after a some period} after the working element fails.
In some cases, such work breaks are unacceptable.
Therefore, the reserve elements may also be in the ``hot'', ``warm'', or ``cold'' reserve.

If reserve elements are in ``cold'' reserve, they do not work, and they are start the work only after the failure of the main element.
This situation coincides with the behaviour of a Queuing system with a finite number of waiting places, and this situation will not be studied here.

If reserve items are in ``hot standby'', this means that they are working at full strength (but idle).
And immediately after the failure of the main element, its work is switched to one of the backup elements (naturally, to the reserve elemeth which at the time of the failure of the main element is in working condition).

It is generally assumed that the distribution of working and recovery time of the backup elements are the same as the corresponding distributions for the primary (main) element.
Working periods or restorable periods of any elements of the system are independent.
Clearly, that such ``hot reservation'' is too expensive.
The behavior of such system is described by the behavior of a Queuing system with a finite number of independently and permanently operating servers; here such system will not be considered.

Finally, in many technical systems, reserve elements are in ``warm reserve'', or partially energized standby redundancy, they work, but not ``in full force''.
I.e. in the case of serviceability of the main element .
In this situation, reliable time of standby element in idle mode is significantly increased, and when the main element fails, the standby element very quickly switches to the working condition.
In this case, generally, the distribution of working time (including being in reserve) and recovery time of reserve element from ``warm reserve'' is different from the corresponding distributions of the main element.

An important characteristic of a restorable element (or system consisting of restorable elements) is {\it reliability factor} $A (t)$, or the probability that at the time $t$ the element (or the system) is in working regime.

For most restorable systems there exists the limit
$$\lim\limits _{t\to \infty}A (t)=A (\infty),$$
and it is very important to know the rate of this convergence.

Note that $A (t)=\mathbb{P} \{X_t\in \mathcal{S}_0\}$, where $X_t$ is the state of the studied system, and $\mathcal{S}_0$ is a subset of the state space of this system.

Therefore, to estimate the rate of convergence of $A (t)$ to its stationary value, it is sufficient to estimate the rate of convergence of the distribution of random variable (r.v.) $X_t$ to its stationary distribution -- in cases where the calculation of $A (t)$ is difficult or impossible.
\section{Restorable element with one warm standby element: case of Markov chain in continuous time}

\subsection{Restorable element without reserve}
For one restorable element with cumulative distribution function (c.d.f.) of working time $F (s)=1-e^{-\lambda s}$, and c.d.f. of recovery time $G (s)=1-e^{-\mu s}$, it is well known that in the case of $A_1 (0)=1$ (see, e.g., \cite[\S 2.3]{GBS})
\begin{equation}\label{A1}
A_{1}(t)\stackrel{\text{\rm def}}{=\!\!\!=}P\{\mbox{at the time $t$ the element works}\}=\frac{\mu + \lambda e^{-(\lambda + \mu)t}}{\mu + \lambda}.
\end{equation}

Accordingly, it is easy to find the convergence rate of availability factor $A_1 (t)$ to its stationary value $A_1(\infty)\stackrel{\text{\rm def}}{=\!\!\!=}\displaystyle\frac{\mu}{\mu + \lambda}$, and this rate is exponential:
\begin{equation}\label{sk}
A_1(t)-A_1(\infty)=\frac {\lambda e^{-(\lambda + \mu)t}}{\mu + \lambda}.
\end{equation}

\subsection{Restorable element with one warm standby element}\label{vossst}

Now consider the restorable element in the situation where working time (and time of the stay in a warm standby state), and recovery time of standby element also have an exponential distributions.

This means that serviceable standby element has an {\it intensity of the failure}, which depends on the state of the main restorable element.
Namely, the first (main) element has c.d.f. of working time $F_1(s)=1-e^{-\lambda_1 s}$, and in the case of its failure, its recovery time has c.d.f. $G_1(s)=1-e^{-\mu_1 s}$; these distributions are independent from the state of reserve element.

That is, if the first (main) element is serviceable at the time $t$, then
\begin{equation}\label{1}
\mathbb{P}\{\mbox{the main element will fail in time $[t, t + \Delta]$}\}=\lambda_1\Delta + o(\Delta).\end{equation}

Similarly, if at time $t$ the main element is faulty, then
\begin{equation}\label{2}
\mathbb{P}\{\mbox{the main element will be serviceable in time $[t, t + \Delta]$}\}=\mu_1\Delta + o(\Delta).
\end{equation}

Constant values of $\lambda_1$ and $\mu_1$ are {\it constant intensities} of failures and recoveries, respectively.

The second (reserve) element has the failure intensity which depends on the state of the main element:
\begin{enumerate}
\item if the main element is serviceable, then the intensity of the failure of reserve element is equal to $\lambda_2$;
\item otherwise, the failure intensity of reserve element is equal to $\Lambda$.
\end{enumerate}
If the reserve element is faulty, then the intensity of its recovery is equal to $\mu_2$ (for simplicity, this intensity is independent of the state of the main restorable element).

It is assumed that at the initial time $t=0$ both elements are serviceable.

The state of the described two-element reliability system at the time $t$ is denoted by the pair $X_t=(i_t,j_t)$.
 Here $i_t,j_t\in \{0;1\}$.
The first component of this pair describes the state of the first element, and the second component of this pair describes the state of the second element.
The serviceable state of each element is indicated by a sign $0$, and the faulty state is indicated by a sign $1$.
The process $X_t$ is a Markov chain in continuous time with the state space $\mathcal{X}=\{(0,0),(1,0), (0,1), (1,1)\}$.

The availability factor of this system is $A_2\stackrel{\text{\rm def}}{=\!\!\!=} \mathbb{P} \{X_t\neq (1,1)\}$.

Denote $p_{i,j}\stackrel{\text{\rm def}}{=\!\!\!=}\mathbb{P}\{X_t= (i,j)\}$.
Analogically \eqref{1}--\eqref{2}, and taking into consideration law of total probability, it implies
$$
\left\{
\begin{array}{l}
p_{0,0}(t + \Delta)=p_{0,0}(1-\lambda_1\Delta - \lambda_2 \Delta) + p_{1,0}\mu_1 \Delta + p_{0,1}\mu_2 \Delta + o(\Delta);
\\
p_{1,0}(t + \Delta)=p_{1,0}(1-\mu_1\Delta - \Lambda\Delta) + p_{0,0}\lambda_1\Delta + p_{1,1}\mu_2 \Delta + o(\Delta);
\\
p_{0,1}(t + \Delta)=p_{0,1}(1-\mu_2\Delta - \lambda_1\Delta) + p_{0,0}\lambda_2\Delta + p_{1,1}\mu_1 \Delta + o(\Delta);
\\
p_{1,1}(t + \Delta)=p_{1,1}(1-\mu_1\Delta -\mu_2\Delta) + p_{1,0}\lambda_2 \Delta + p_{0,1}\lambda_1 \Delta + o(\Delta).
\end{array}
\right.
$$
Considering, that $p_{0,0} + p_{1,0} + p_{0,1} + p_{1,1}=1$, it holds the Kolmogorow equations:
\begin{equation}\label{eq}
\left\{
\begin{array}{l}
\dot p_{0,0}(t)=p_{0,0}(-\lambda_1 - \lambda_2 ) + p_{1,0}\mu_1 + p_{0,1}\mu_2 ;
\\
\dot p_{1,0}(t)=p_{0,0}(\lambda_1-\mu_2) + p_{1,0}(-\mu_1 - \Lambda-\mu_2) -p_{0,1}\mu_2 + \mu_2 ;
\\
\dot p_{0,1}(t)=p_{0,0}(\lambda_2 -\mu_1) -p_{1,0} \mu_1 + p_{0,1}(-\mu_2 - \lambda_1-\mu_1) + \mu_1.
\end{array}
\right.
\end{equation}
As $p_{0,0}(t) + p_{0,1}(t)=A_1(t)$ is the availability factor of the main element, $\alpha_1\stackrel{\text{\rm def}}{=\!\!\!=}-(\lambda_1+\mu_1)$ is the eigenvalue of the system \eqref{eq} -- see \eqref{A1}.

The other two eigenvalues are solutions of the square equation:
\begin{multline}\label{korni}
\alpha_{2,3}=\displaystyle\frac{-\Lambda}{2}- \displaystyle\frac{\lambda_1}{2} - \displaystyle\frac{\lambda_2 }{2}-
\displaystyle\frac{\mu_1}{2}-\mu_2\pm
\\ \\
\pm \displaystyle\frac {\sqrt{(\Lambda + \mu_1)^2 + (\lambda_1 + \lambda_2 )^2 + 2\mu_1(\lambda_1-\lambda_2 + 2\mu_2-8)- 2\Lambda(\lambda_1 + \lambda_2 )}}{2}.
\end{multline}
Note, that in this case the convergence rate of the availability factor $A_2(t)=p_{0,0}(t) + p_{0,1}(t) + p_{1,0}(t)$ to the stationary value
$$
A_2(\infty)={ \textstyle\frac{ \mu_1(\mu_1\lambda_2+\Lambda(\lambda_1+\lambda_2 )+2(\mu_2^2+\mu_1-\mu_2)) + \mu_2(\lambda_1(\lambda_1+\lambda_2 +\mu_2) +\mu_1\mu_2(2\lambda_1+\lambda_2 +\Lambda) )} {(\lambda_1 + \mu_1) (\mu_1(2 + \lambda_2 ) + \mu_2(\mu_2 + \Lambda + \lambda_1+\lambda_2 ) + \Lambda(\lambda_1 + \lambda_2))}},
$$
can be less than in (\ref{sk}).

Really,
$$
A_2(t)= A_2(\infty)+ C_1 e^{\alpha_1 t}+ C_2 e^{\alpha_2 t}+ C_3 e^{\alpha_3 t},
$$
where $C_1$, $C_2$, $C_3$ depend on the initial state of the system.
If $\Lambda$, $\lambda_2$ and $\mu_2$ small enough, then according to \eqref{korni},
$$
\alpha_{2}\stackrel{\text{\rm def}}{=\!\!\!=}\displaystyle\frac{-\Lambda}{2}- \displaystyle\frac{\lambda_1}{2} - \displaystyle\frac{\lambda_2 }{2}-
\displaystyle\frac{\mu_1}{2}-\mu_2+ \displaystyle\frac {\sqrt{\;\;\ast\;\;\ast\;\;\ast\;\;}}{2}>\alpha_1,
$$
and the rate of convergence of $ A_2 (t)$ to the limit value is not less than $ C_2 e^{\alpha_2 t}$ -- of course, in the case where $C_2\neq 0$.

It should be reminded that an important problem of the reliability theory is to estimate the rate of convergence of the availability factor for the stationary value not only in the simplest (exponential) case described above, but also in more complicated situations.

In this paper, we consider the case where the distributions of the working times and recovery times of both (main and reserve) elements can depend on each other, and their distribution is described by {\it intensities bounded above and below}.

To estimate the rate of convergence of the availability factor to the stationary value, the rate of convergence of the distribution of the {\it full} state of this restorable system of in the total variation metrics will be estimated.
The estimation of this convergence rate automatically implies the estimation of the convergence rate of the availability factor.

To obtain an estimate of the convergence rate of the distribution of the {\it full} state of the restorable system in the metric of the total variation, the {\it coupling method} and the markovization of the semi-Markov process described the behavior of the restorable system by means of some extantion of the state space of this system will be used.
Recall that the calculations will be done under the assumption that the distributions of working times and recovery times of both elements are described by {\it bounded above and below intensities.}

In this case, it will be shown that the convergence rate of the parameters of the restorable system is exponential.
The results of this kind have many important applications -- see., e.g., \cite{kal1,kal2}.

\section{Restorable element with one warm standby element: semi-Markov case}
So, consider a restorable system consisting of two elements.

It is assumed that the distributions of working times and recovery times may be not exponential, but they are absolutely continuous.

Now, {\it the full state of the system} is two pairs:
\begin{equation}\label{proc}
X_t\stackrel{\text{\rm def}}{=\!\!\!=}((i_t,x_t),(j_t,y_t)),
\end{equation}
where $i_t=0$ or $i_t=1$ if the first (main) element is serviceable or not at the time $t$ correspondingly.
And the value $x_t$ is equal to the time elapsed from the last change in the state of $i_t$ of the first (main) element to the time $t$.

The pair $ (j_t, y_t)$ describes the state of the reserve element at the time $t$ in the same way.

Note that the behavior of the restorable system with reserve element(s) is usually described by a semi-Markov (piecewise linear) process (see, e.g. \cite{GK}).

However, using the {\it intensities} instead of their distribution functions, allows to consider this process as a Markov process $X_t$ on the state space $\mathcal{X}\stackrel{\text{\rm def}}{=\!\!\!=}\{0,1\} \times R_+\times \{0,1\}\times R_+$.

\subsection{Description of process behavior by {\it intensities}}
A significant number of publications on Queuing theory and reliability theory use the {\it constant intensity} of input flow of the customers, of the service completion or of failures.

The constant intensity corresponds to an exponential distribution (time between the inputs of the customers, or of service times, or of working times, or of recovery times). The behavior of the models in the Queueing theory and in the reliability theory with {\it the constant intensities} can be described by the Markov chain in continuous time, as it was shown in the \eqref{vossst} section.
However, in many applied problems it is necessary to consider arbitrary distributions of random periods, and in this case the concept of {\it intensity of completion of a random period of the time} is used (see, e.g. \cite{GBS, GK}).

\begin{definition}\label{def1}
Let $\mathcal{T}$ be some random period of time (with c.d.f. $\Phi(x)$), started at the time $t=0$.
And let this period $\mathcal{T}$ not end at the time $s>0$.
Then
$$\mathbb{P} \{\mathcal{T} \in
(s,s+\Delta]| \mathcal{T}>s\}=\displaystyle
\frac{\Phi(s+\Delta)-\Phi(s)}{1-\Phi(s)}.
$$

If c.d.f. $\Phi(x)$ is absolutely continuous, then
$$
\displaystyle\phi(s)\stackrel{\text{\rm def}}{=\!=} \lim\limits
_{\Delta \downarrow 0}\frac{\mathbb{P} \{\mathcal{T} \in
(s,s+\Delta]| \mathcal{T}>s\}}{\Delta} =
\frac{\Phi'(s+0)}{1-\Phi(s)}.
$$

The function $\phi(s)$ is called { \texttt{the intensity}} of the end
of the random period $\mathcal{T}$ at the time $s$ given
$\mathcal{T}\geqslant s$. \hfill \ensuremath{\triangleright}
\end{definition}

\begin{remark}\label{-intens} 
The definition \ref{def1} implies:
$$
\mathbb{P} \{\mathcal{T}\in(s,s+\Delta]| \mathcal{T}>s\}=\phi(s)\Delta+ o(\Delta).\eqno{\triangleright}
$$
\end{remark}
\begin{remark}\label{zv9} 
C.d.f. and probability density function the random variable $\mathcal{T}$ can be calculated by the intensity of the end
of the random period $\mathcal{T}$:
\begin{equation}\label{zv2} 
\Phi(s)=1-\exp\left(-\int\limits_0^s \phi(u)\,\mathrm{d}\, u\right); \qquad \Phi'(s)= \phi(s)\exp\left(-\int\limits_0^s \phi(u)\,\mathrm{d}\, u\right).
\end{equation}

It is easy to see that: if c.d.f. the random variable $\mathcal{T}$ is exponential, then the intensity of the end of this period is constant, and Vice versa.

If c.d.f. $\Phi(s)$ is not continuous, then it can also consider the intensity of the end of the period with c.d.f. $\Phi(s)$, but the formulas (\ref{zv2}) take an other form, taking into account the jumps of c.d.f. \hfill \ensuremath{\triangleright}
\end{remark}

\begin{definition}
The random variable $\eta$ does not exceed the random variable $\theta$ \texttt{by distribution} if the inequality
$$
\displaystyle F_\eta(s)=\mathbb{P} \{\eta\leqslant s\}
\geqslant \mathbb{P} \{\theta\leqslant s\}=F_\theta(s)
$$
is true for all $s\in\mathbb{R}$.

In other words, c.d.f. of the random variable $\theta$ does not exceed c.d.f. of the random variable $\eta$.
It is the order in the set of the random variables.
Denote it by $\eta\prec \theta$ -- see \cite{sht}.\hfill \ensuremath{\triangleright}
\end{definition}
\begin{remark}\label{rem3}
If the intensity $\phi(s)$ of the end of the random period $\mathcal{T}$ satisfies the inequality
$$
0<c<\phi(s)<C<\infty,
$$
then
\begin{equation}\label{sravni}
\mathcal{T}_-\prec \mathcal{T}\prec\mathcal{T}_+,
\end{equation}
where $\mathbb{P}\{\mathcal{T}_-\leqslant s\}=1-e^{-C s}$, and $\mathbb{P}\{\mathcal{T}_-\leqslant s\}=1-e^{-c s}$.
Correspondingly, $\mathbb{E}\,\mathcal{T}\leqslant \displaystyle\frac{1}{c}$.

Indeed, in accordance with \eqref{zv2}, $\mathbb{P}\{\mathcal{T}_-\leqslant s\}=1-\exp\left(-\displaystyle\int\limits_0^s C\,\mathrm{d}\,u\right)$, and $\mathbb{P}\{\mathcal{T}_-\leqslant s\}=1-\exp\left(-\displaystyle\int\limits_0^s c\,\mathrm{d}\,u\right)$, whence it follows \eqref{sravni}.
\hfill \ensuremath{\triangleright}
\end{remark}
Finally, after this preparatory information, we describe the behavior of the considering system in the language of intensities.

\subsection{Assumptions and denotations}
Assume that the first (main) element has a failure intensity (in the operating mode) and a recovery intensity (in the non-operating mode) {\it depending on the complete state of the system} -- it is $\lambda_1(X_t)$ and $\mu_1 (X_t)$ correspondingly.

Similarly, the second (reserve) element has a failure intensity $ \lambda_2 (X_t)$, and the recovery intensity of the second (reserve) element in case of failure is $\mu_2 (X_t)$.

Assume that when the first (main) element is restored, it immediately starts to work, and the reserve element, even if it was serviceable, goes into a state of warm reserve.

In other words,
\begin{equation}\label{prx}
\left\{
\begin{array}{l}
\mathbb{P}\{i_{t+\Delta}=1, x_{t+\Delta}\in[0,\Delta), j_{t+\Delta}=j_t, y_{t+\Delta}=y_t+\Delta | i_t=0\}=
\\
\hspace{7cm}=\lambda_1(X_t)\Delta+0(\Delta);
\\
\mathbb{P}\{i_{t+\Delta}=0, x_{t+\Delta}\in[0,\Delta), j_{t+\Delta}=j_t, y_{t+\Delta}=y_t+\Delta | i_t=1\}=
\\
\hspace{7cm}=\mu_1(X_t)\Delta+0(\Delta);
\\
\mathbb{P}\{j_{t+\Delta}=0, y_{t+\Delta}\in[0,\Delta), i_{t+\Delta}=i_t, x_{t+\Delta}=x_t+\Delta | j_t=0\}=
\\
\hspace{7cm}=\lambda(X_t)\Delta+0(\Delta);
\\
\mathbb{P}\{j_{t+\Delta}=1, y_{t+\Delta}\in[0,\Delta), i_{t+\Delta}=i_t, x_{t+\Delta}=x_t+\Delta | j_t=1\}=
\\
\hspace{7cm}=\mu_2(X_t)\Delta+0(\Delta).
\end{array}
\right.
\end{equation}

The equalities (\ref{prx}) specify the transition function of the process $X_t$, i.e. the $X_t$ process is a Markov process with a state space $\mathcal{X}$ and with transition probabilities \eqref{prx}.

The model described here includes the previously considered exponential case, as well as the situation where in the case of a long duration the work of the main element, the reserve element begins to work in advance in a more intensive mode, in order to quickly enter into operation in the event of a failure of the main element.
Also, the work of the main element may be less intense if the reserve element is faulty, to exclude the possibility of simultaneous failure of both elements, etc.

\begin{center}
{\bf Assumptions. }
\end{center}

{\it Suppose, that}
\begin{equation}\label{usl}
\begin{array}{c}
0<\lambda_m^-\leqslant \lambda_m(X_t)\leqslant {\lambda_m^+}<\infty;
\\
\\
0<{\mu_m^-}\leqslant \mu_m(X_t)\leqslant {\mu_m^+}<\infty.
\end{array}
\end{equation}
\hfill \ensuremath{\triangleright}

Let $\xi_m^+$, $\xi_m^-$, $\zeta_m^+$, $\zeta_m^-$ be an additional random variables with c.d.f.:
\begin{equation}\label{oboz}
\begin{array}{ll}
\mathbb{P} \{\xi_m^+\leqslant s\}=1-e^{\int\limits _0^s -\lambda_m^+ (u) \,\mathrm{d} u}; & \mathbb{P} \{\xi_m^-\leqslant s\}=1-e^{\int\limits _0^s -\lambda_m^- (u) \,\mathrm{d} u};
\\
\mathbb{P} \{\zeta_m^+\leqslant s\}=1-e^{\int\limits _0^s -\mu_m^+ (u) \,\mathrm{d} u}; & \mathbb{P} \{\zeta_m^-\leqslant s\}=1-e^{\int\limits _0^s -\mu_m^- (u) \,\mathrm{d} u}.
\end{array}
\end{equation}
In the future, the sets of such random variables $\xi_m^+(k)$, $\xi_m^-(k)$, $\zeta_m^+(k)$, $\zeta_m^-(k)$, $k\in \mathbb{N}$ will be used, assuming that these random variables are mutually independent.
\begin{remark}\label{rem4}
Denote: $\xi_m(k)$ is $k$-th working time of $m$-th element ($m=1,2$), and $\zeta_m(k)$ is $k$-th recovery time of $m$-th element.

According to Remark \ref{rem3},
$$
\zeta_m(k)\prec \zeta_m^-, \qquad \xi_m(k)\prec \xi_m^-. \eqno{\triangleright}
$$
\end{remark}

\subsection{Ergodicity of process $X_t$}
Denote $\mathcal{P}_t$ -- distribution of process $X_t$ at the time $t$, i.e. $\mathcal{P}_t(A)=\mathbb{P}\{X_t\in A\}$, $A\in \mathcal{B}(\mathcal{X})$.
Clearly, $\mathcal{P}_t$ depends on $X_0$.

The process $X_t$ is not {\it regenerative}, and standard methods of the proof of its ergodicity as an ergodicity of the regenerative process are not applicable (see \cite{zv1,zv2}).

\begin{lemma}\label{lemma}
Under conditions \eqref{usl}, the process $X_t$ is ergodic, i.e. there exists an invariant probability measure $\mathcal{P}$ on the state space $\mathcal{X}$ -- such that $\mathcal{P}_t\Longrightarrow\mathcal{P}$ for all initial states $X_0$ of the process $X_t$.
\end{lemma}
\textbf{Proof. }
Denote (for $m=1,2,\ldots$)
$$
\begin{array}{cc}
\theta_1 \stackrel{\text{\rm def}}{=\!=} \inf\{t>0:\; i_t=0\}, & \theta_1' \stackrel{\text{\rm def}}{=\!=} \inf\{t>\theta_1:\; i_t=1\}, \\
\theta_{m+1} \stackrel{\text{\rm def}}{=\!=} \inf\{t>\theta_m':\; i_t=0\}, & \theta_{m+1}' \stackrel{\text{\rm def}}{=\!=} \inf\{t>\theta_{m+1}':\; i_t=0\}.
\end{array}
$$

Here $\theta_k$ is the time of the $k$-th repair of the main element.

Remarks \ref{rem3} and \ref{rem4} imply:
\begin{equation}\label{sr}
\mathbb{E}\,\theta_1\leqslant \displaystyle\frac{1}{\lambda^-_1}+\displaystyle\frac{1}{\mu^-_1};\qquad \mathbb{E}\,(\theta_{k+1}-\theta_k)\leqslant \displaystyle\frac{1}{\lambda^-_1}+ \displaystyle\frac{1}{\mu^-_1}.
\end{equation}

Let $\varepsilon>0$ be small enough.
Denote
$$\mathcal{S}^\varepsilon\stackrel{\text{\rm def}}{=\!\!\!=} \{(1,x,1,y):\, \max\{x,y\}<\varepsilon\}; \qquad\mathcal{S}^\varepsilon\subset \mathcal{X}.
$$

At the time $\theta_k$ ($k>1$) the first (main) restorable element hits into working state; at this time, $i_{\theta_k}=0$, $x_{\theta_k}=0$.

Two cases are possible with respect to the second (reserve) element at the time $\theta_k$.

1. If $j_{\theta_k}=1$, i.e. the second element is not serviceable, then the residual time of its stay in not-serviceable state {\it by distribution} is less than the random variable $\zeta_2^-$,
So, the second element will hits in a serviceable state for a time less than $\varepsilon$  with the probability greater then $\varpi_1 \stackrel{\text{\rm def}}{=\!=} 1-e^{-\varepsilon\mu^-_2}$. Thus, before the time $\theta_k+\varepsilon$, the process $X_t$ will hit into the set $\mathcal{S}^\varepsilon$.

2. If $j_{\theta_k}=0$, i.e. the second element is serviceable, then the residual time of its stay in serviceable state is less then the random variable $\xi^-$ {\it by distribution}, and the subsequent recovery time less then the random variable $\zeta_2^-$ by {\it distribution}.
Thus,
\begin{multline*}
\mathbb{P}\left\{
\begin{array}{c}
                   \mbox{the second element will fail faster than in the time $\frac\varepsilon 2$,}
                   \\
                   \mbox {and then recover faster than in the time $\frac\varepsilon2$}
                 \end{array}
\right\}
\geqslant
\\
\\
\geqslant\varpi_2 \stackrel{\text{\rm def}}{=\!=} \left(1-e^{-\frac{\varepsilon}{2}\mu^-_2} \right)\left(1-e^{-\frac{\varepsilon}{2}\lambda^-_2 }\right).
\end{multline*}
So, the process $X_{t}$ will be in the small subset $\mathcal{S}^\varepsilon$ of the state space $\mathcal{X}$ before the time $\vartheta_k+\varepsilon$ with probability greater then $\max\{\varpi_1 ,\varpi_2 \}$.

Hence, the time between consequent hits of the $X_t$ process to the set $\mathcal{S}^\varepsilon$ is the geometric sum of the random variable with finite means (see \eqref{sr}), and it has a finite expectation.
Therefore, the ergodicity of the process $X_t$ follows from the principle of Harris-Khasminsky.
\hfill \ensuremath{\blacksquare}
\section{Main result}
\begin{theorem}\label{th}
If conditions \eqref{usl} are satisfied, then for any initial state $X_0$ of the process $X_t$, {\it it is possible to calculat} the numbers $\alpha>0$ and $\mathfrak{K}=\mathfrak{K}(\alpha, \lambda^-_m, \lambda^+_m,\linebreak \mu^-_m,\mu^+_m)$ such that
$
\|\mathcal{P}_t-\mathcal{P}\|_{TV}\leqslant \mathfrak{K}e^{-\alpha t}
$
for all $t\geqslant 0$.

Here $\|\cdot\|_{TV}$ is the total variation metrics.
\end{theorem}
Recall that the availability factor $A_2(t)$ is
\begin{equation}\label{af}
A_2(t)=\mathbb{P}\{X_t\notin\{((1,\cdot),(1,\cdot))\}\}= 1-\mathcal{P}_t((1,\cdot),(1,\cdot)).
\end{equation}

\begin{corollary}
From \eqref{af} it follows that in the conditions of the Theorem \ref{th} the inequality
$$|A_2(t)-A_2(\infty)|\leqslant \mathfrak{K}e^{-\alpha t}$$
is true. \hfill \ensuremath{\triangleright}
\end{corollary}

\noindent{\bf Proof of the Theorem \ref{th}.} This proof is the algorithm of computing the values $\alpha>0$ and $\mathfrak{K}$. It  is based on the {\it coupling method} for Markov processes -- see, e.g. \cite{zv6,zv14}.

Here, we only show  the main ideas of this proof.

Recall that the ``classical'' coupling method is not applicable to our model, since its behavior is described by the Markov process \emph{in continuous time}.
Instead, the concept of successful coupling proposed in \cite{zv5} will be used.

The concept of ``{\it  the common part of distributions}'' will be also used.
The common part of two distributions with the densities $f_1(\cdot)$, $f_2(\cdot)$ is a number $\kappa \stackrel{\text{\rm def}}{=\!=} \displaystyle\int\limits_\mathbb{R}\min\{f_1(u), f_2(u)\}\,\mathrm{d}u$.
Also,  the {\it Basic Coupling Lemma} in the simplest form will be used (see, e.g., \cite{zv13,verbut}).
About  application of the coupling method in Queuing theory, see, e.g., \cite{zv7,GZ}).

\subparagraph{1. Successful coupling.}
Let $X_t$ and $\widehat{X}_t$ be two independent Markov processes with the same transition function \eqref{prx}, but with different initial states at time $t=0$.

Suppose that (dependent) processes $Y_t=((i_t,x_t),(j_t,y_t))$ and $\widehat{Y}_t=((\widehat{\imath}_t,\widehat{x}_t),(\widehat{\jmath}_t,\widehat{y}_t))$ are constructed on a certain probability space, in such a way that:

\begin{enumerate}
\item $Y_t\stackrel{\mathcal{D}}{=}X_t$ and $\widehat{Y}_t\stackrel{\mathcal{D}}{=}\widehat{X}_t$ for all \emph{non-random} $t$;
\item $\mathbb{P} \{\tau(X_0,\widehat{X}_0)<\infty\}=1$, where $\tau(X_0,\widehat{X}_0)=\tau(Y_0,\widehat{Y}_0)=\inf\{t>0:\; Y_t=\widehat{Y}_t\}$.
\end{enumerate}

Recall that this pair of  processes $Y_t=((i_t,x_t), (j_t,y_t)$ and $\widehat{Y}_t=((\widehat{\imath}_t,\widehat{x}_t),(\widehat{\jmath}_t,\widehat{y}_t))$ is called \emph{successful coupling} for the processes $X_t$ and $\widehat{Y}_t$, and $\tau (X_0,\widehat{X}_0)$ is called \emph{ coupling epoch}.

Now, our goal is a construction of the successful coupling and an estimation of exponential moments of a random variable $\tau(X_0,\widehat{X}_0)$.

\subparagraph{2. Construction of processes $Y_t$ and $\widehat{Y}_t$.}
In ordered to construct the processes $Y_t$ and $\widehat{Y}_t$ let us consider a probability space $\prod\limits_{k=0}^\infty (\Omega_k,\mathbb{P} _k,\mathcal{B}_k)$, where $(\Omega_k,\mathbb{P} _k,\mathcal{B}_k)$ are some probability spaces which my be specified.
In the probability space $(\Omega_k,\mathbb{P} _k,\mathcal{B}_k)$, the required for successful coupling random variables will be defined.

We will proceed ``step by step'' subsequent stopping times $t_1$, $t_2$, $\ldots$, where one of the components $i,j,\widehat{\imath},\widehat{\jmath}$ of the pair $(Y_t,\widehat{Y}_t)$ changes.

Let the state of both processes be known at the moment $t_k$. Then, we know the intensities of the change of the components $i_t,j_t,\widehat{\imath}_t,\widehat{\jmath}_t$ at the time $t_k$, and the values $i_{t_k},j_{t_k},\widehat{\imath}_{t_k},\widehat{\jmath}_{t_k}$ at the moment $t_k$.
So, joint probability distributions of the residual times of the stay of the elements at their states $i,j,\widehat{\imath},\widehat{\jmath}$ are known.

At this step, one of the probability spaces $(\Omega_k,\mathbb{P} _k,\mathcal{B}_k)$ is used.
On this probability space, we construct our four random variables in such a way that their joint probability distribution coincides with the described above distribution.
Now, let us define the minimal value of these four random variables. 
Denote it by $\varsigma_k$.

At the moment $t_{k+1}\stackrel{\text{\rm def}}{=\!\!\!=} t_k+\varsigma_k$, the next change of one of the components $i,j,\widehat{\imath},\widehat{\jmath}$ occurs.

The above procedure is now repeated at the time $t_{k+1}$, by induction.

Note that the proposed procedure of construction of stochastic processes $Y_t$ and $\widehat{Y}_t$ does not guarantee their coincidence  at any time.

\subparagraph{3. Basic Coupling Lemma (see, e.g., \cite{zv13,verbut}).}
Here the simplest formulation of the Basic Coupling Lemma is given (see, e.g., \cite{GZ1}).

\begin{lemma}\label{osn}
If the random variable $\vartheta_1$ and $\vartheta_2$ have c.d.f. $\Phi_1(s)$ and $\Phi_2(s)$ cor\-res\-pon\-din\-g\-ly, and their common part $\kappa\stackrel{\text{\rm def}}{=\!\!\!=} \int\limits_\mathbb{R} \min\{\Phi_1'(s), \Phi_2'(s)\}\,\mathrm{d}\, s>0$, then it can construct (on some probability space) the random variables $\widehat{\vartheta}_1$ and $\widehat{\vartheta}_2$ such, that
\begin{enumerate}
\item $\widehat{\vartheta}_1\stackrel{\mathcal{D}}{=} {\vartheta}_1$, $\widehat{\vartheta}_2\stackrel{\mathcal{D}}{=} {\vartheta}_2$;
\item $\mathbb{P} \{\widehat{\vartheta}_1= \widehat{\vartheta}_2\}=\kappa$. \hfill \ensuremath{\triangleright}
\end{enumerate}
\end{lemma}
\begin{remark}\label{lem+}
The statement of Lemma \ref{osn} is naturally transferred to any finite number of random variables.

Let $\vartheta_1$, $\vartheta_2$, $\ldots$, $\vartheta_n$ be the random variable with probability dencities $\varphi_1(s)$, $\varphi_2(s)$, $\ldots$, $\varphi_n(s)$ correspondingly, and  $\kappa \stackrel{\text{\rm def}}{=\!\!\!=}\int\limits_\mathbb{R}\min\{\varphi_1(s), \varphi_2(s), \ldots, \varphi_n(s)\}>0$.
Then on some probabilistic space it is possible to construct the random variables $\widehat{\vartheta}_1(s)$, $\widehat{\vartheta}_2(s)$, $\ldots$, $\widehat{\vartheta}_n(s)$ such that
\begin{enumerate}
\item $\widehat{\vartheta}_i\stackrel{\mathcal{D}}{=} {\vartheta}_i$, $i=1,2,\ldots n$;
\item $\mathbb{P} \{\widehat{\vartheta}_1= \widehat{\vartheta}_2=\ldots=\widehat{\vartheta}_n\}=\kappa$.

\noindent{\bf Proof} of this fact is repeated by the proof of Lemma \ref{osn} from \cite{GZ1}.\hfill \ensuremath{\triangleright}
\end{enumerate}
\end{remark}

\subparagraph{4. Coupling of the processes $Y_t$ and $\widehat{Y}_t$.}

Firstly, consider the times $\theta_k$.
 $\theta_k$ are the times of recovery of the main element of the process $Y_t$; $\theta_k\prec \sum\limits_{i=1}^k(\xi_1^-(i)+\zeta_1^-(i))$ (see Lemma \ref{lemma}).

Similar to the proof of the Lemma \ref{lemma}, for some fixed $\varepsilon>0$ we can affirm the following:

\begin{enumerate}
\item With a probability greater than
$$\pi_1\stackrel{\text{\rm def}}{=\!\!\!=} (1-e^{-\frac\varepsilon2\lambda_2^-})(1-e^{-\frac\varepsilon2\mu_2^-}),
$$
the process $Y_t$ hits to the set $\mathcal{S}^\varepsilon$ in the time interval $[\theta_k,\theta_k+\varepsilon]$.
\item With a probability greater than
$$\pi_2\stackrel{\text{\rm def}}{=\!\!\!=} ( 1-e^{-\frac\varepsilon2 \lambda_1^-})( 1-e^{-\frac\varepsilon2 \mu_1^-})( 1-e^{-\frac\varepsilon2 \Lambda^-})( 1-e^{-\frac\varepsilon2 \lambda_2^-})( 1-e^{-\frac\varepsilon2 \mu_2^-})
$$
the process $\widehat{Y}_t$ hits to the set $\mathcal{S}^\varepsilon$ in the time interval $[\theta_k,\theta_k+\varepsilon]$.
\end{enumerate}
Hence, with a probability greater than $\pi_1\pi_2$, both processes ${Y}_t$ and $\widehat{Y}_t$ will hit to the set $\mathcal{S}^\varepsilon$ in the time interval $[\theta_k,\theta_k+\varepsilon]$.

In this set $\mathcal{S}^\varepsilon$ all elements of both system are serviceable, and the residual time of the stay of all elements in serviceable condition has a common part
\begin{multline*}
\kappa_1\stackrel{\text{\rm def}}{=\!\!\!=}\int \limits_0^\infty \min \{\lambda_1(X_{\theta_k+\varepsilon+s}) e^{-\lambda_1(X_{\theta_k+\varepsilon+s})}, \lambda_2(X_{\theta_k+\varepsilon+s})e^{-\lambda_2(X_{\theta_k+\varepsilon+s})} \}\,\mathrm{d}\, s \geqslant
\\
\geqslant\frac{\min\{\lambda_1^-,\lambda_2^-\}} {\max\{\lambda_1^+,\lambda_2^+\}}>0.
\end{multline*}

Therefore,  it is possible to construct four random variables (the residual time of the stay of all elements in serviceable condition) by such a way, that they will coincide with probability no less than $\kappa_1$, on  some probability space.
That is, both processes will be simultaneously in the state $(0,0,0,0)$ through time $\eta$, which is less than the random variable with with c.d.f. $\Phi(s)=1-e^{-\min\{\lambda_1^-,\lambda_2^-,\mu_1^-,\mu_2^-\}s}$.

In fact, it can first to construct the coincidence of the states of the first (main) elements of both processes, and then it can  to construct the coincidence the states of the reserve elements, but the description of this procedure is much more complex.

For the constructin of the coupling, it can also consider the times $\theta_k'$ -- the times of failures of the first (main) element.

Again, repeating the reasoning of the Lemma \ref{lemma}, it can estimate the probability $\mathbb{P}\{\mbox{on the intervale $[\theta_k',\theta_k'+\varepsilon]$ all elements of both systems will be faulty}\}$, and then it can apply the modified Basic Coupling Lemma \ref{osn}, and show that with probability greater than some $\kappa_2>0$, the processes will coincide.

The proposed here estimates of the probabilities of the entering into certain set for the pair of processes are not optimal.
The calculations can be improved by considering different combinations of the states of both elements of both processes and different ways of choosing the value of $\varepsilon$ and other part of this value instead of $\frac{\varepsilon}{2}$.

\subparagraph{Conclusion.} So, with some probability $\tilde{\kappa}$ (with a very rough estimate), it can be argued that on any interval $[\theta_k,\theta_{k+1}]$, the processes $Y_t$ and $\widehat{Y}_t$ will coincide.
Moment of coincidence is $\tau=\tau(X_0,\widehat{X}_0)$.

Therefore, it can apply the basic coupling inequality for the constructed successful coupling:

\begin{multline}\label{zv12} 
|\mathbb{P} \{X_t\in S\} - \mathbb{P} \{\widehat{X}_t\in S\}|= |\mathbb{P} \{Y_t\in S\} - \mathbb{P} \{\widehat{Y}_t\in S\}|=
\\
=|\mathbb{P} \{Y_t\in S\;\&\; \tau>t\} - \mathbb{P} \{\widehat{Y}_t\in S\;\&\; \tau>t\} +
\\
+ |\mathbb{P} \{Y_t\in S\;\&\; \tau<t\} - \mathbb{P} \{\widehat{Y}_t\in S\;\&\; \tau<t\}|=
\\
=|\mathbb{P} \{Y_t\in S\;\&\; \tau>t\} - \mathbb{P} \{\widehat{Y}_t\in S\;\&\; \tau>t\} |\leqslant \mathbb{P} \{\tau>t\}\leqslant
\\
\leqslant \mathbb{P} \{e^{\alpha\tau}>e^{\alpha t}\}\leqslant \frac{\mathbb{E} \,e^{\alpha\tau}}{e^{\alpha t}},
\end{multline}
as for $\tau<t $ ``coupling'' of the processes $Y_t$ and $\widehat{Y}_t$ already happened, and their distributions at the time $t$ are equal, and $\mathbb{P} \{Y_t\in S \,\&\, \tau<t\}=\mathbb{P} \{\widehat{Y}_t\in S\,\&\, \tau<t\}$; here $S\in\mathcal{B} (\mathcal{X} )$.

Thus, it remains to estimate $\mathbb{E} \, e^{\beta\tau}$ and find out at what $\beta>0$ this expectation is finite.

Considering that there is a uniform (independent of initial conditions) estimate
$$
\tau(X_0,\widehat{X}_0)\prec\theta_1+\sum\limits _{k=1}^\nu (\theta_{k+1}-\theta_k)+\eta\prec\eta+\sum\limits_{k=1}^\nu(\xi_k^-+\zeta_k^-),
$$
where $\theta_{k+1}-\theta_k\prec \xi_1^-+\zeta_1^- \stackrel{\mathcal{D}}{=}
\xi_k^-+\zeta_k^-$, and also $\mathbb{P} \{\nu>n\}\leqslant (1-\tilde{\kappa}) ^n$, it can estimate the value of $\alpha<\tilde{\kappa}$ and $\mathfrak{K}=\mathbb{E} \,e^{\alpha\tau}$.

As inequality \eqref{zv12} is uniforme by the initial states of the processes $X_t$ and $\widehat{X}_t$, this inequality is true in the situation where the process $\widehat{X}_t$ has an initial distribution equal to the stationary distribution $\mathcal{P}$.
Therefore,
$$
\sup\limits_{S\in \mathcal{B}(\mathcal{X})}|\mathbb{P} \{X_t\in S\}-\mathcal{P}(S)|=\|\mathcal{P}_t(S)-\mathcal{P}(S)\|_{TV}\leqslant \frac{\mathbb{E} \,e^{\alpha\tau}}{e^{\alpha t}}=\mathfrak{K}e^{-\alpha t}.
$$

\section{Conclusion}
The constuction of successful coupling presented here is significantly simplified and the bounds, which can be obtained according to the described scheme, are rather rough.

However, using the proposed approach with the analysis of all possible situations suitable for possible constructins of coupling epoch, as well as the possibility to take into account some specific characteristics of the considered system, and using of different initial states of the process $X_t$ allow to significantly improve the proposed estimate.

\paragraph{Acknowledgement.} The author thanks E.Yu.~Kalimulina, V. V.~Kozlov and A.Yu. Ve\-re\-ten\-nikov for valuable recommendations in the preparation of the article.
The work is supported by RFBR (project No 17-01-00633 A).